\documentclass[final,5p,times,twocolumn]{elsarticle}

\usepackage{amsmath} 
\usepackage{csvsimple} 
\usepackage{siunitx}   
\usepackage{booktabs}  
\usepackage{afterpage}

\journal{Swarm and Evolutionary Computation}

\begin{document}

\begin{frontmatter}



\title{QUASAR: \\ An Evolutionary Algorithm \\ to Accelerate High-Dimensional Numerical Optimization}


\author{Julian G. Soltes} 

\affiliation{organization={Regis University},
            addressline={3333 Regis Blvd.}, 
            city={Denver},
            postcode={80221}, 
            state={Colorado},
            country={USA}}

\begin{abstract}
High-dimensional numerical optimization presents a persistent challenge in computational science. This paper introduces Quasi-Adaptive Search with Asymptotic Reinitialization (QUASAR), an evolutionary algorithm to accelerate convergence in complex, non-differentiable problems afflicted by the curse of dimensionality.

QUASAR expands upon the core principles of Differential Evolution (DE), introducing quasi-adaptive mechanisms to dynamically balance exploration and exploitation in its search. Inspired by the behavior of quantum particles, the algorithm utilizes three highly stochastic mechanisms that augment standard DE: 1) probabilistic mutation strategies and scaling factors; 2) rank-based crossover rates; 3) asymptotically decaying covariance reinitializations.

Evaluated on the notoriously difficult CEC2017 benchmark suite of 29 test functions, QUASAR achieved the lowest overall rank sum (367) using the Friedman test, outperforming DE (735) and L-SHADE (452). Geometric mean comparisons show average final solution quality improvements of $3.85 \times$ and $2.07 \times$ compared to DE and L-SHADE, respectively ($p \ll 0.001$), with average optimization speed averaging $1.40 \times$ and $5.16 \times$ faster. QUASAR's performance establishes it as an effective, efficient, and user-friendly evolutionary algorithm for complex high-dimensional problems.
\end{abstract}



\end{frontmatter}



\section{Introduction}
High-dimensional numerical optimization remains a persistent challenge across many fields, including machine learning, science, engineering, and finance \cite{zhou2024evolutionary}. This challenge stems from the curse of dimensionality; the exponential growth of the search space renders gradient-based and parametric search methods ineffective \cite{shan2008survey}. Consequently, optimization metaheuristics, particularly evolutionary algorithms, have become essential for locating optimal solutions in complex function landscapes.

Among these methods, Differential Evolution (DE) has established itself as a robust and highly effective optimizer due to its simplicity, efficiency, and effective mutation strategies \cite{storn1996usage}. However, standard DE variants can struggle to maintain an effective balance between global exploration and local exploitation, particularly in high-dimensional search spaces. This often leads to stagnation or premature convergence that requires frequent restarts or extensive hyperparameter tuning, a significant practical drawback for real-world applications.

To address these limitations, this paper proposes Quasi-Adaptive Search with Asymptotic Reinitialization (QUASAR), an evolutionary algorithm designed to dynamically balance the exploration-exploitation trade-off and achieve superior performance in high-dimensional spaces. This method is conceptually inspired by the search for stability in a stellar core, where quantum particles probabilistically explore configurations searching for the most stable state.

QUASAR extends the DE framework by incorporating a highly stochastic multi-strategy approach, including probabilistically selected mutation strategies and factors, as well as an adaptive crossover rate based on solution fitness. The algorithm's primary innovation is an asymptotically decaying reinitialization mechanism, utilizing the spatial covariance matrix of top-performing solutions to inject high-quality genetic diversity into the population.

The efficacy of QUASAR is demonstrated through rigorous evaluation against standard DE and its state-of-the-art variant L-SHADE across the full CEC2017 test function benchmark suite. The results show large and statistically significant improvements in final solution quality, with high computational efficiency. These performance gains are seen using its default parameters, highlighting the ease of use and scalability of the algorithm.

The remainder of the paper is structured as follows: Section 2 details the QUASAR algorithm, including its mutation, crossover, selection, and reinitialization mechanisms. Section 3 defines the experimental design and statistical analysis, Section 4 analyzes the results. Section 5 describes the ease of use and implementation of the algorithm. Section 6 summarizes the results and concludes the paper, with Section 7 proposing future enhancements and validations.

\section{QUASAR Algorithm}
Quasi-Adaptive Search with Asymptotic Reinitialization (QUASAR) operates using a highly stochastic, multi-strategy approach. The mechanisms defined below synergize to create a dynamic 'push-and-pull' effect on the solution vectors.

\subsection{Initial Population}
The initial population is generated as a Sobol sequence by default, due to its low discrepancy and superior uniformity in high-dimensional spaces compared to other quasi-Monte Carlo sampling methods \cite{bonneel2025sobol}. The default population size, $N$, is set to $10\times D$, where $D$ is the number of dimensions.

Additional initialization methods include uniform Latin Hypercube and random distributions, as well as non-uniform Hyperellipsoid Density sequences \cite{hyperellipsoid}.

\subsection{Mutation}
Solutions are recombined with one another, analogous to particles sharing information via quantum entanglement. This 'spooky action at a distance' (A. Einstein, 1947) serves as the conceptual basis for the mutation mechanisms below.

QUASAR uses three distinct mutation strategies, each being a modified Differential Evolution variant. These strategies are selected dynamically based on the specified entanglement probability, \textit{entangle rate}. This value defaults to $0.33$, resulting in the three strategies being applied equally across the population.

\subsubsection{Local Mutation (Exploitation):}
The local mutation strategy, '\textit{Spooky-Best}', is chosen with probability \textit{entangle rate}. This strategy is similar to the DE/best/1 format, exploiting the immediate search space around the current best solution $\mathbf{X}_{\text{best},g}$.
\begin{equation}
    \mathbf{v}_{i,g} = \mathbf{X}_{\text{best},g} + F_{\text{local}}(\mathbf{X}_{i,g} - \mathbf{X}_{\text{rand},g})
\tag{\textit{Spooky-Best}}
\end{equation}

The associated mutation factor $F_{\text{local}}$ is sampled from a normal distribution $\mathcal{N}(0, 0.33^2)$, allowing for small bi-directional perturbations about the best solution.

\subsubsection{Global Mutations (Exploration):}
Two global mutation strategies are applied to the remaining $(1-\textit{entangle rate})$ fraction of the population. Promoting global exploration, these strategies, '\textit{Spooky-Current}' and '\textit{Spooky-Random}', each have a 50\% probability of being selected:
\begin{gather}
    \mathbf{v}_{i,g} = \mathbf{X}_{i,g} + F_{\text{global}}(\mathbf{X}_{\text{best},g} - \mathbf{X}_{\text{rand},g}) \tag{\textit{Spooky-Current}} \\
    \mathbf{v}_{i,g} = \mathbf{X}_{\text{rand},g} + F_{\text{global}}(\mathbf{X}_{i,g} - \mathbf{X}_{\text{rand},g}) \tag{\textit{Spooky-Random}}
\end{gather}

The mutation factor $F_{\text{global}}$ is sampled from a bimodal distribution $\mathcal{N}(0.5, 0.25^2) + \mathcal{N}(-0.5, 0.25^2)$ to drive non-local exploration. It is worth clarifying that the \textit{Spooky-Random} strategy uses the same random vector $\mathbf{X}_{\text{rand}}$ twice.
\\

The terms used in the three mutation strategies are defined as:
\begin{itemize}
    \item $\mathbf{v}_{i,g}$: Mutant vector for the $i$-th individual at generation $g$.
    \item $\mathbf{X}_{i,g}$: Current solution vector.
    \item $\mathbf{X}_{\text{best},g}$: Best known solution in the population.
    \item $\mathbf{X}_{\text{rand},g}$: Randomly chosen solution vector.
    \item $F_{\text{local}}$ and $F_{\text{global}}$: Corresponding mutation factors.
\end{itemize}

\begin{figure}[!ht]
\centering
\includegraphics[width=0.45\textwidth]{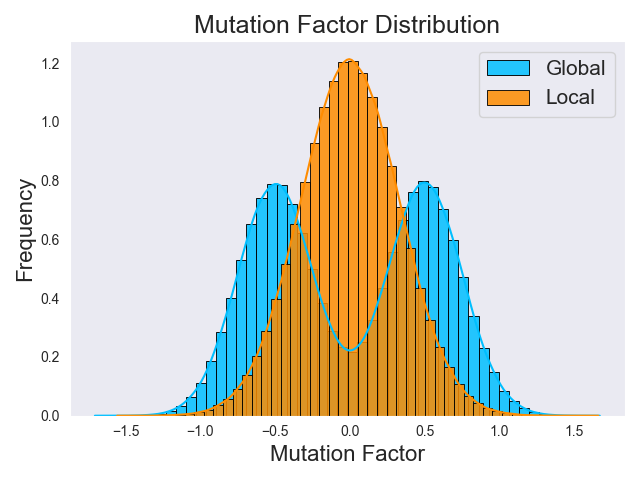}
\caption{Mutation factor distributions for global and local mutation strategies.}
\label{fig:mutations}
\end{figure}

\subsection{Crossover}
QUASAR uses binomial crossover with a dynamic crossover rate ($CR_i$) inversely proportional to the solution's fitness rank. The best solutions are rapidly recombined to maximize both exploration and exploitation through synergy with global and local mutations.

The raw crossover rate ($CR_\text{raw},i$) for the $i$-th solution is calculated based on its rank ($\textit{rank}_i \in [0,N-1]$), where $\textit{rank}=0$ corresponds to the best solution (lowest objective function value). The final dynamic crossover rate $CR_i$ is capped at a minimum of $0.33$. This mechanism ensures that the best solutions have a higher chance of inheriting components from the new mutant vector $\mathbf{v}_i$.

\begin{equation}
CR_{\text{raw},i} = \frac{N - 1 - \text{rank}_i}{N - 1}
\end{equation}

\begin{equation}
CR_i = \max(CR_{\text{raw},i}, 0.33)    
\end{equation}

The trial vector $\mathbf{u}_i$ is then generated component-wise ($n$ dimension) using binomial crossover, comparing a uniform random number $\text{rand}(0,1)$ against $CR_i$:
\begin{equation}
    u_{i,n} = \begin{cases} v_{i,n}, & \text{if } \text{rand}(0,1) \leq CR_i \\ X_{i,n}, & \text{if } \text{rand}(0,1) > CR_i \end{cases}
\end{equation}

\subsection{Selection}
A greedy elitism strategy is used to determine the next generation's population, where the target vector for generation $g+1$ is the better of the two solutions between the current vector $\mathbf{X}_{i,g}$ and the trial vector $\mathbf{u}_{i,g}$. If a solution undergoes reinitialization, it bypasses the selection step.
\begin{equation}
    \mathbf{X}_{i,(g+1)} = \begin{cases} \mathbf{u}_i, & \text{if } f(\mathbf{u}_i) < f(\mathbf{X}_i) \\ \mathbf{X}_i, & \text{if } f(\mathbf{u}_i) \geq f(\mathbf{X}_i) \end{cases}
\end{equation}

\subsection{Asymptotic Reinitialization}
The covariance-guided vector reinitialization applies to a fixed 33\% of the worst-performing solutions. In early generations, these solutions have a high probability of 'tunneling' to more promising regions. This mechanism complements the rank-based crossover; while the best solutions recombine frequently, lower quality solutions explore locally until stochastically tunneling to a better location.

The probability of reinitialization $P_{\text{reinit}}(g)$ decays asymptotically from $P_{\text{reinit}}(0)=1$ toward a target final probability, calculated using the current generation $g$:

\begin{equation}
    P_{\text{reinit}}(g) = \exp\left(\frac{\ln(P_{\text{final}})}{g_{\text{final}} \cdot g_{\text{max}}} \cdot g\right)
\end{equation}

where $g_\text{max}$ is the maximum generation count, $P_{\text{final}}=0.33$ is the target final probability, and $g_\text{final}=0.33$ defines the fraction of generations at which $P_{reinit}$ reaches $P_{final}$. This ensures the probability drops to $P_{\text{final}}$ at $g=g_\text{max} \cdot g_{\text{final}}$.

During a reinitialization event, new solutions are sampled from a multivariate distribution shaped by top performing individuals. The top $M$ solutions (where $M=25\% \text{of the population}$) are used to calculate the mean vector $\boldsymbol{\mu}$ and the covariance matrix $\boldsymbol{\Sigma}$:

\begin{equation}
    \boldsymbol{\mu} = \frac{1}{M} \sum_{i=1}^{M} \mathbf{X}_i
\end{equation}

\begin{equation}
\boldsymbol{\Sigma} = \frac{1}{M-1} \sum_{i=1}^{M} (\mathbf{X}_i - \boldsymbol{\mu})(\mathbf{X}_i - \boldsymbol{\mu})^{\top} + \epsilon \mathbf{I}
\end{equation}

where $\epsilon$ is a small scalar ($10^{-12}$) and $\mathbf{I}$ is the identity matrix, added to prevent $\boldsymbol{\Sigma}$ from being singular. Similarly, the asymptotic reinitialization is bypassed for cases where $N<D$. The effects of skipping the reinitialization can be seen in the $N=100$ sample-variant trials (Fig. \ref{fig:gmerf_samples}).

New vectors $\mathbf{y}$ are generated by sampling from a multivariate normal distribution, with noise added to further promote exploration.

\begin{equation}
    \mathbf{y} \sim \mathcal{N}(\boldsymbol{\mu}, \boldsymbol{\Sigma}) +\boldsymbol{\delta}
\end{equation}

The noise vector $\boldsymbol{\delta}$ is sampled from a multivariate normal distribution, with respect to the search bounds:

\begin{equation}
\boldsymbol{\delta} \sim \mathcal{N}\left(0, \left(\frac{\mathbf{b}_{\text{high},n} - \mathbf{b}_{\text{low},n}}{20}\right)^2 \mathbf{I}\right)
\end{equation}

Finally, all components of the new solution $\mathbf{X}_{\text{new}}$ are clipped to ensure they remain within the defined search bounds $\mathbf{B} =[\mathbf{b}_{\text{low}}, \mathbf{b}_{\text{high}}]$.

\begin{equation}
    \mathbf{X}_{\text{new},n} = \text{clip}(\mathbf{y}_n,\mathbf{B}_n)
\end{equation}

These new solutions completely replace the reinitialized solutions and are not subjected to crossover mechanisms.

Due to the covariance calculations only using a fraction of the population, along with decaying probability of occurrence, the additional overhead becomes negligible compared to the objective function evaluations. For extreme cases, the reinitialization can be disabled, which is the default behavior for $N<D$ problems. 

\begin{figure}[h!]
\vspace{-0.33cm}
\centering
\includegraphics[width=0.45\textwidth]{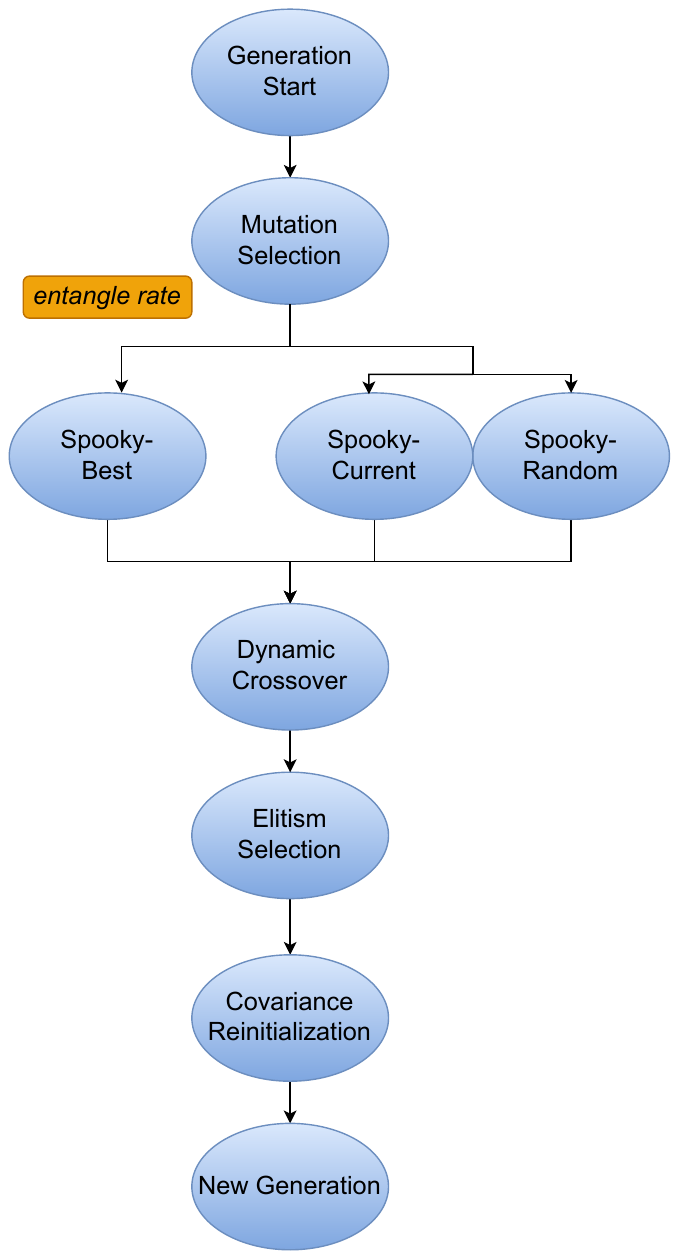}
\vspace{-1.5cm}
\caption{QUASAR workflow for evolving each generation.}
\label{fig:workflow}
\end{figure}

\section{Experimental Design}
\subsection{Optimizers \& Test Functions}
The performance of QUASAR was benchmarked against modern implementations of Differential Evolution (DE) and L-SHADE, specifically SciPy's DE and MealPy's L-SHADE. 

This choice is justified by the fact that the QUASAR algorithm is built on the DE framework, as well as DE's standing as a robust and effective optimizer \cite{price1996differential}. L-SHADE was specifically chosen due to its success in CEC optimization competitions  \cite{mo2019efficient}.

The trials were conducted using all 29 functions of the CEC2017 numerical optimization benchmark suite.

\subsection{Experimental Setup}
To ensure a fair comparison, each optimizer was run using their default hyperparameters. The initial population generation methods were also left to defaults; the tested sample sizes not being powers of 2 diminishes the uniformity of QUASAR's initial Sobol sequence \cite{joe2008constructing}, providing a theoretical initial handicap.

Since SciPy's DE implementation internally scales the input population size by dimensionality ($N \times D$  when all dimensions have equal bounds), a correction factor of $1/D$ is applied to its initial sample size.

The computational budget for metaheuristic comparison is often defined by the maximum number of objective function evaluations ($FES_{max}$) \cite{kazikova2021objective}. Due to L-SHADE’s linear population reduction, the $FES_{max}$ becomes a difficult comparison metric to compare against. As such, the primary comparisons are performed on final solution fitness errors and overall optimization times, given the same number of generations, using modern implementations via SciPy’s DE and MealPy’s L-SHADE. 

By default, the L-SHADE implementation features tolerance parameters for early convergence. Additionally, the QUASAR and DE implementations include a final polishing step, using Powell and L-BFGS-B minimization, respectively. These were all disabled for the trials to ensure the optimizers fully utilize their computational budget, terminating only when reaching the maximum number of generations $g_{max}$.

\subsection{Experimental Trials:}
The trials were structured into two distinct phases to assess QUASAR’s behavior under varying dimensionalities and population sizes. 

The results are compared based on the overall aggregated rank sums, geometric mean error factors ($R_{\text{error}}$), and optimization run time ratios $R_{\text{time}}$.

\subsubsection{Dimension-Variant Trials:}
The dimension-variant trials assessed QUASAR’s scalability by scaling the problem dimension $D$ with a fixed population size ($N=1000$). $30$ independent trials were performed for each test function/dimension combination to ensure statistical significance. 

The dimensions were varied between $[10,30,50,100]$, resulting in a total of 10,400 full optimization runs; 3,480 for each algorithm. The computational budget was kept constant, with maximum number of iterations $g_{max}=100$ and fixed sample size $N=1000$.

\subsubsection{Sample-Variant Trials:}
The sample-variant trials assessed the impact of population size on performance. The computational budget is varied between population sizes of $[100, 250, 500, 1000, 5000]$. $30$ independent trials were performed for each test function / sample size combination, using maximum number of iterations $g_{max}=100$ and fixed dimensionality $D=100$. 

The $N=5000$ trials, however, were only run for 10 trials to ensure statistical significance while mitigating the high computational overhead, resulting in a total of 11,310 full optimization runs; 3770 for each optimization algorithm.

\begin{table}[h!]
\centering
\caption{Experimental Parameters for Sample and Dimension Variations}
\label{tab:test_parameters}
\begin{tabular}{l|cc}
\toprule
 & \textbf{Sample-variant} & \textbf{Dimension-variant} \\
\midrule \midrule
\textbf{D}      & 100                          & [10, 30, 50, 100] \\
\textbf{N} & [100, 250, 500, 1000, 5000] & 1000 \\
$\mathbf{g_{\max}}$ & 100                  & 100 \\
\bottomrule
\end{tabular}
\end{table}

\subsection{Statistical Analysis:}
\subsubsection{Overall Rank:}
To draw the overall statistical comparison between the optimizers across all function/dimension and function/sample size combinations, the Friedman test was used. 

For each unique combination, the median final error value $f$ across all independent trials was calculated for each optimizer. The optimizers were ranked based on these median errors, with the best performance (lowest error) receiving a rank of 1.

The ranks were summed across all scenarios to produce the aggregated rank sum for each optimizer. A lower sum indicates superior overall performance. 

\subsubsection{Solution Quality Improvements:}
The primary comparison metric used to quantify the magnitude of performance improvement is the geometric mean error reduction factor $R_{\text{error}}$. This provides an unbiased measure of relative final objective penalty $f$ for the comparison optimizer to QUASAR across all $T$ independent trials for each scenario:

$$R_{\text{error}} = \sqrt[T]{ \prod_{i=1}^{T} \frac{f_{\text{comp},i}}{f_{\text{QUASAR}, i}}}$$

Here, $f_{\text{comp},i}$ and $f_{\text{QUASAR},i}$ are the final error values achieved by the competitor algorithm (e.g., DE or L-SHADE) and QUASAR, respectively, in the $i$-th trial.

For the overall results, the average $\bar{R}_{\text{error}}$ across all data is calculated as the geometric mean of the $R_{\text{error}}$ obtained for each of $k$ function/dimension or function/sample size combinations.

$${\bar{R}_{\text{error}}} = \sqrt[k]{ \prod_{j=1}^{k}{R_{\text{error},j}}}$$

\subsubsection{Computational Efficiency:} 
The computational efficiency of the optimizers is quantified by comparing their average optimization run times. 

The run times ratios $R_{\text{time}}$ for QUASAR vs. DE \& L-SHADE are shown in Table \ref{tab:runtime_data}, by dimension and population size. $R_{\text{time}}$ is calculated by dividing the arithmetic average run time $\bar{t}$ of the comparison method by the average run time of $\text{QUASAR}$ across all trials for that specific dimension.

$$R_{\text{time}} = \frac{\bar{t}_{\text{comp}}}{\bar{t}_{\text{QUASAR}}} = \frac{\sum_{i=1}^{T} t_{\text{comp}, i}}{\sum_{i=1}^{T} t_{\text{QUASAR}, i}}$$

The overall efficiency metric, $\bar{R}_{\text{time}}$, is defined as the geometric mean of the $R_{\text{time}}$ ratios across all $k$ function/dimension or function/sample size combinations.

\begin{equation} \bar{R}_{\text{time}} = \sqrt[k]{ \prod{j=1}^{k} R_{\text{time}, j} } \end{equation}

\subsubsection{Statistical Significance:}
 For direct comparison within each dimension or sample size, the paired Wilcoxon signed-rank test is performed. This non-parametric test was chosen as the final error distributions are on different scales for each function, are assumed not to be normally distributed, and to ensure consistency with the non-parametric ranking methodology. 
 
 The signed-rank test was conducted on the final error values ($f$) and run times to generate the \textit{p}-values for each function/dimension combination.

\section{Results}
The primary solution quality metrics are the Friedman aggregated rank sum and the geometric mean error reduction factor $R_{\text{error}}$, quantifying QUASAR's improvements in final solution fitness. The computational efficiency is measured by comparing the ratio of full optimization run times $R_{\text{time}}$.

QUASAR demonstrated top performance in the rank sum test, shown in Table \ref{tab:rank_sums}. Similarly, it demonstrated consistent improvements in geometric mean error reduction factor $R_\text{error}$ (Table \ref{tab:gmerf_data}) and run times $R_\text{time}$ (Table \ref{tab:runtime_data}). 

\begin{table}[h]
    \centering
    \caption{Rank Sums by Sample Size and Dimension}
    \label{tab:rank_sums}
    \begin{tabular}{l|ccc}
        \toprule
         & \textbf{N} & \textbf{D} & \textbf{Overall} \\
        \midrule \midrule
        \textbf{QUASAR}  & 217 & 150 & 367 \\
        \textbf{L-SHADE} & 223 & 229 & 452 \\
        \textbf{DE}      & 430 & 305 & 735 \\
        \bottomrule
    \end{tabular}
\end{table}

\begin{table*}[h!]
    \centering
    \caption{Overall Performance Comparison (Geometric Mean Ratios)}
    \label{tab:overall}
    \begin{tabular}{l|ccc|ccc}
        \toprule
         & \multicolumn{3}{c|}{\textbf{vs L-SHADE}} & \multicolumn{3}{c}{\textbf{vs DE}} \\
         & $\mathbf{\bar{R}}$ & \textbf{95\% CI} & \textbf{\textit{p}-value} & $\mathbf{\bar{R}}$ & \textbf{95\% CI} & \textbf{\textit{p}-value} \\
        \midrule \midrule
        \textbf{Quality} ($\mathbf{\bar{R}_{\text{error}}}$) & 2.07 & (0.62, 6.96) & $5.30 \times 10^{-144}$ & 3.85 & (0.93, 16.03) & $5.78 \times 10^{-136}$ \\
        \textbf{Run Time} ($\mathbf{\bar{R}_{\text{time}}}$)  & 5.16 & (4.57, 5.83) & $3.36 \times 10^{-141}$ & 1.40 & (1.36, 1.44)  & $4.52 \times 10^{-23}$ \\
        \bottomrule
        \multicolumn{7}{l}{\footnotesize \textit{Ratios > 1 indicate superior QUASAR performance.}}
    \end{tabular}
\end{table*}

\subsection{Solution Quality}
The ranks and geometric mean fitness errors $R_{\text{error}}$ of the final solutions are calculated for both the dimension-variant and sample-variant trials. Solution quality comparisons are shown in Table \ref{tab:gmerf_data}.

\subsubsection{Dimension-Variant:}
The dimension-variant analysis assessed QUASAR's convergence against DE and L-SHADE across four discrete dimensions $[10,30, 50, 100]$, using a fixed sample size $N=1000$ and maximum generations $g_{max}=100$.

\begin{figure}[h!]
\centering
\includegraphics[width=0.45\textwidth]{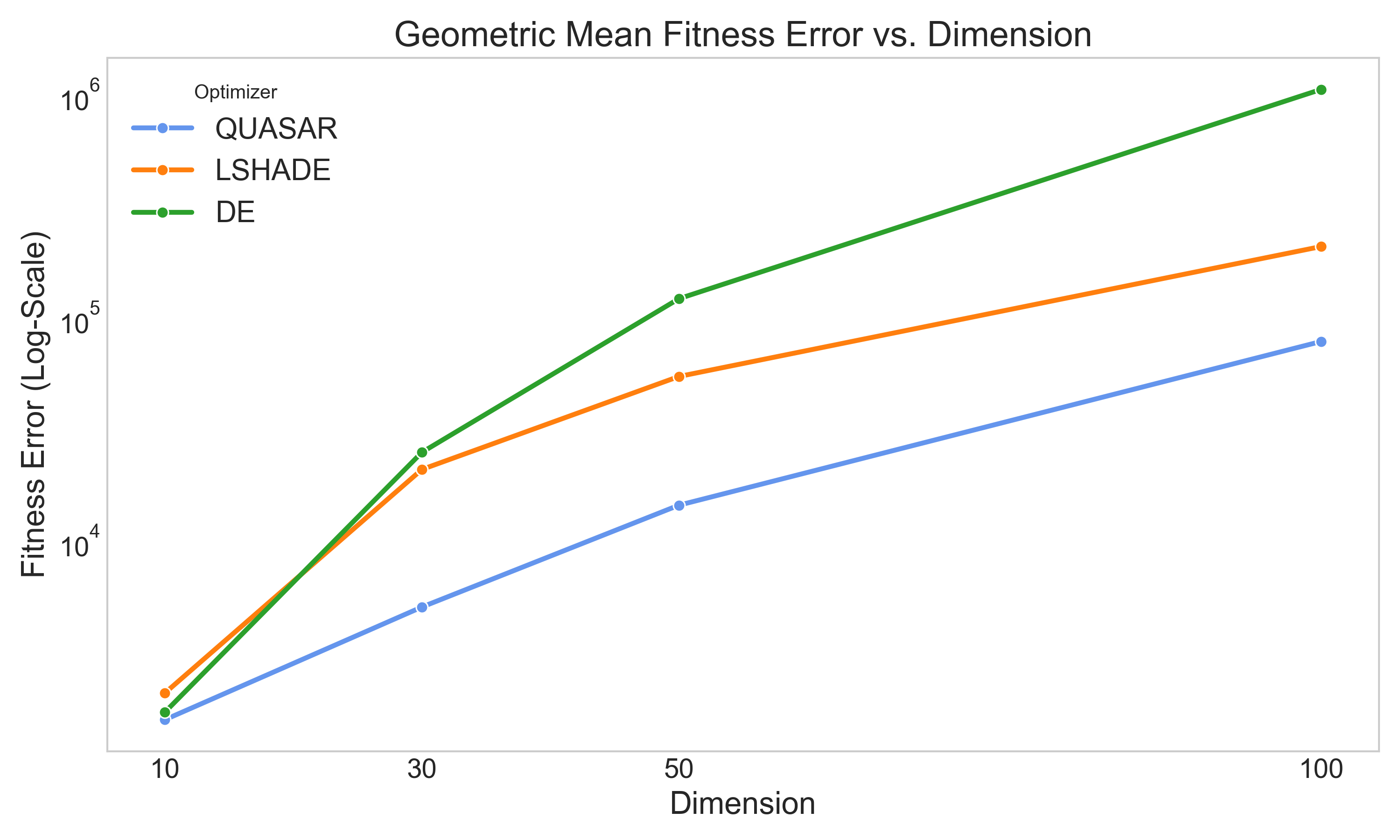}
\caption{Geometric mean error (log-scale) by dimension for each optimizer, with population size of 1000.}
\label{fig:gmerf_dim}
\end{figure}

The geometric mean error reduction factors ($R_{\text{error}}$) were consistently high across all dimensions (Fig. \ref{fig:gmerf_dim}). The 95\% confidence intervals validate its consistent improvements; the intervals do not fall below $R_{\text{error}} = 1.0 \times$. The largest improvements were seen against DE in $100D$, achieving a $R_{\text{error}}$ of $\text{13.52} \times$ with a CI upper bound of $16.68\times$. The comparison against L-SHADE showed the highest improvements in $30D$, with $R_{\text{error}} = 4.14$, reaching CI upper bounds of 4.74.

\subsubsection{Sample-Variant:}
To assess the impact of varying population sizes, the optimizers were compared across a fixed $100D$ and 100 maximum generations $g_{max}$, using various sample sizes ($N$).

\begin{figure}[!ht]
\centering
\includegraphics[width=0.45\textwidth]{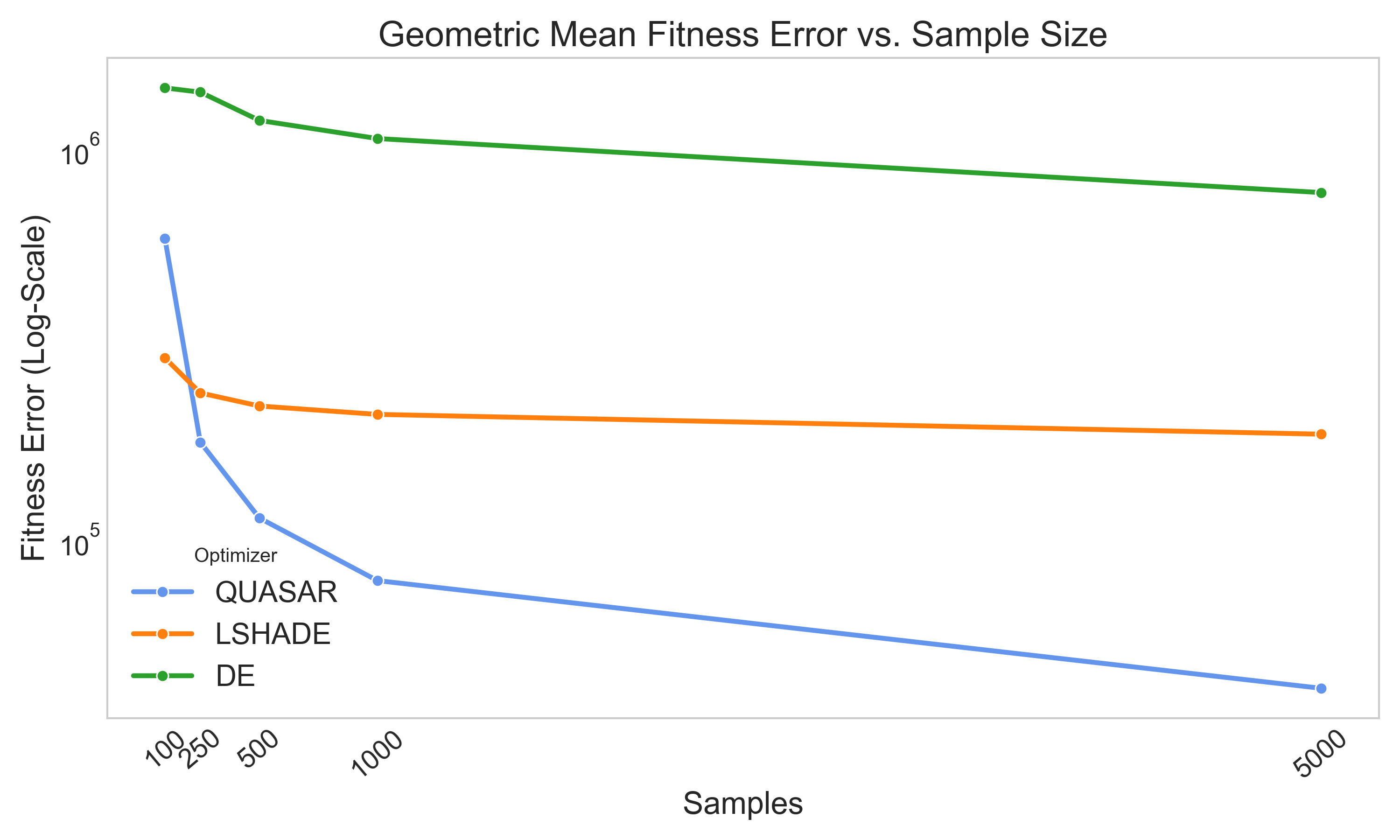}
\caption{Geometric mean error (log-scale) by population size for each optimizer, in 100 dimensions.}
\label{fig:gmerf_samples}
\end{figure}

The fitness improvement factors saw the highest gains with higher population sizes; in the $N=5000$ case, QUASAR reached a geometric mean error improvement factor $R_{\text{error}}$ of $\text{18.5} \times$ compared to DE, and $4.47  \times$ compared to L-SHADE. These performance increases scale with population size (Fig. \ref{fig:gmerf_samples}). The only scenario where QUASAR did not outperform L-SHADE was the $N=100$ case, where L-SHADE achieved $R_{\text{error}} = 0.49 \times$, highlighting the impact of bypassing the covariance reinitialization for $N<D$.

\begin{table*}[h!]
    \centering
    \caption{Solution Quality Comparison (Geometric Mean Error)}
    \label{tab:gmerf_data}
    \begin{tabular}{c|ccc|ccc}
        \toprule
         & \multicolumn{3}{c|}{\textbf{vs L-SHADE}} & \multicolumn{3}{c}{\textbf{vs DE}} \\
        & $\mathbf{R_{error}}$ & \textbf{95\% CI} & \textbf{\textit{p}-value} & $\mathbf{R_{error}}$ & \textbf{95\% CI} & \textbf{\textit{p}-value} \\
        \midrule \midrule
        $\mathbf{D=10}$  & 1.30 & (1.20, 1.40) & $7.0 \times 10^{-16}$  & 1.06 & (1.01, 1.12) & $7.0 \times 10^{-3}$    \\
        $\mathbf{D=30}$  & 4.14 & (3.62, 4.74) & $5.0 \times 10^{-91}$  & 4.94 & (4.24, 5.77) & $3.0 \times 10^{-108}$  \\
        $\mathbf{D=50}$  & 3.78 & (3.32, 4.31) & $9.0 \times 10^{-79}$  & 8.44 & (7.09, 10.04) & $5.0 \times 10^{-134}$ \\
        $\mathbf{D=100}$ & 2.67 & (2.32, 3.06) & $2.0 \times 10^{-36}$  & 13.52 & (10.95, 16.68) & $2.0 \times 10^{-141}$ \\
        \midrule
        $\mathbf{N=100}$  & 0.49 & (0.44, 0.56) & $7.0 \times 10^{-86}$  & 2.43 & (2.04, 2.89) & $3.0 \times 10^{-46}$   \\
        $\mathbf{N=250}$  & 1.34 & (1.19, 1.50) & $3.3 \times 10^{-2}$   & 7.86 & (6.51, 9.49) & $3.0 \times 10^{-141}$  \\
        $\mathbf{N=500}$  & 1.93 & (1.70, 2.20) & $2.0 \times 10^{-8}$   & 10.40 & (8.54, 12.67) & $2.0 \times 10^{-140}$ \\
        $\mathbf{N=1000}$ & 2.66 & (2.32, 3.06) & $5.0 \times 10^{-34}$  & 13.52 & (10.95, 16.68) & $2.0 \times 10^{-141}$ \\
        $\mathbf{N=5000}$ & 4.47 & (3.41, 5.86) & $3.0 \times 10^{-37}$  & 18.54 & (12.49, 27.53) & $3.0 \times 10^{-48}$  \\
        \bottomrule
    \end{tabular}
\end{table*}

\subsection{Computational Efficiency}
The computational efficiency of QUASAR was measured using the overall ratio of run time $R_{\text{time}}$ (seconds) relative to DE and L-SHADE. The full set of comparisons is shown in Table \ref{tab:runtime_data}.

\subsubsection{Dimension-Variant Efficiency:}
\begin{figure}[!ht]
\centering
\includegraphics[width=0.45\textwidth]{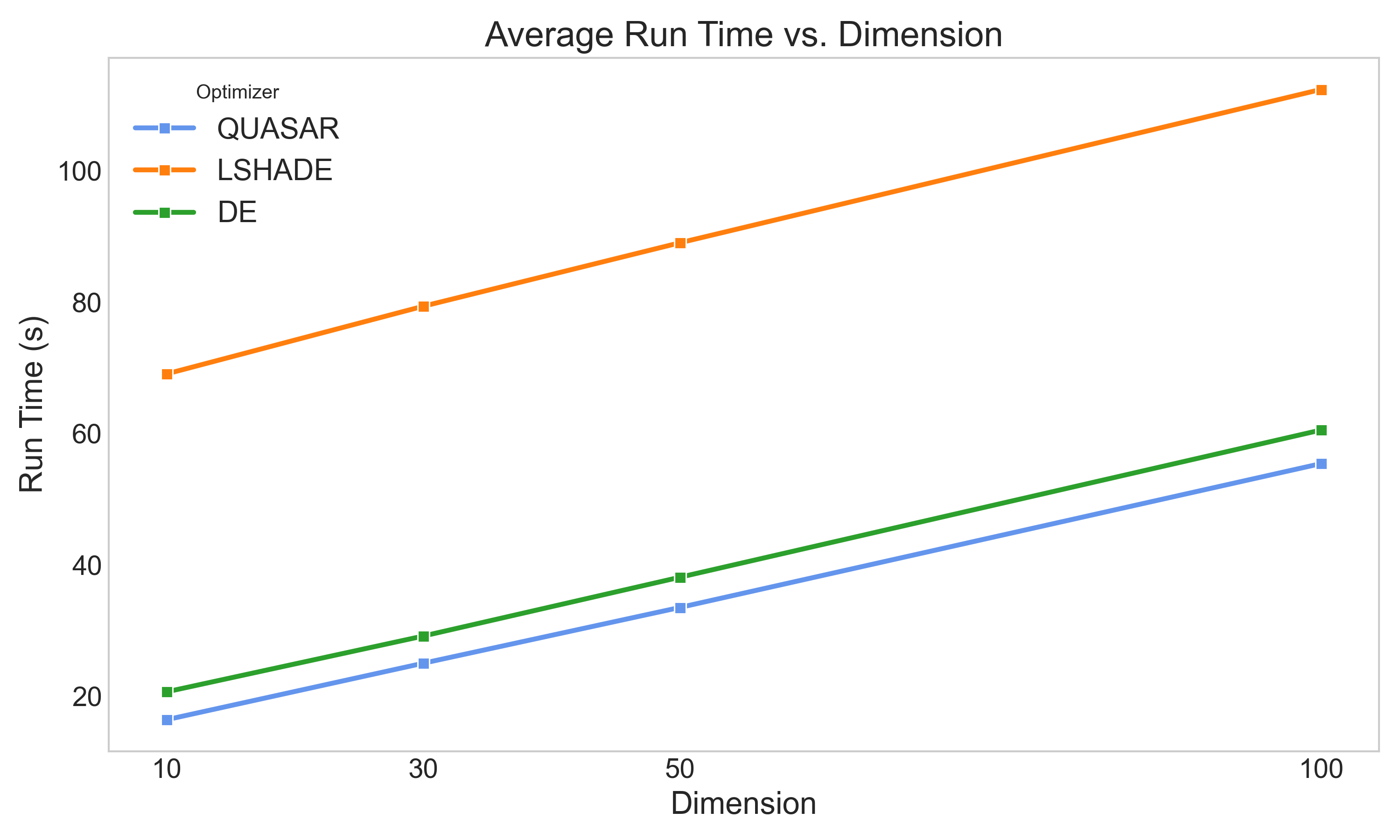}
\caption{Average optimization run times for each optimizer, by dimension.}
\label{fig:runtime_dimensions}
\end{figure}

In the dimension-variant trials, QUASAR demonstrated consistently faster run times. The results, shown in Fig. \ref{fig:runtime_dimensions},  consistently show that QUASAR achieved faster optimization speeds across all dimensions tested, scaling linearly with dimensionality. 

\subsubsection{Sample-Variant Efficiency:}
\begin{figure}[!ht]
\centering
\includegraphics[width=0.45\textwidth]{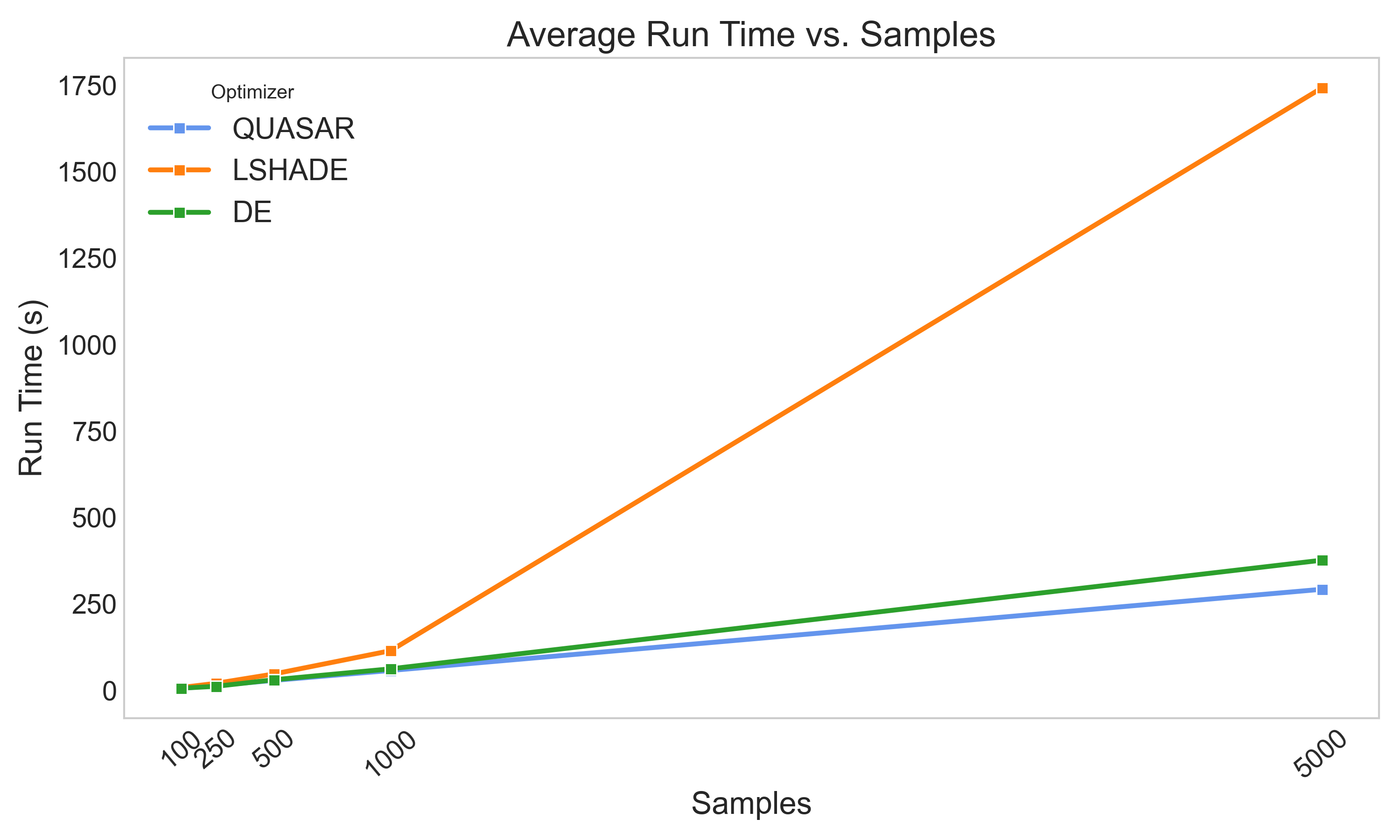}
\caption{Average optimization run times for each optimizer, by population size.}
\label{fig:runtime_samples}
\end{figure}

In the sample-variant trials, run times were consistently faster than DE and L-SHADE, averaging $7.78 \times$ and $1.39 \times$ speed increases respectively across all sample-variant trials (Fig. \ref{fig:runtime_samples}).

The $N=250$ trials, however, demonstrated a ratio $R_{\text{time}}$ vs. DE of $0.84 \times$, indicating that DE was approximately 16\% faster. This likely represents the overhead of the covariance matrix calculation relative to the low total evaluation count, whereas the covariance calculations are bypassed in the lower $N=100$ case where QUASAR was faster than DE.

\begin{table*}[t!]
    \centering
    \caption{Efficiency Comparison (Geometric Mean Run Time)}
    \label{tab:runtime_data}
    \begin{tabular}{c|cc|cc|ccc}
        \toprule
         & \multicolumn{2}{c|}{\textbf{vs L-SHADE}} & \multicolumn{2}{c|}{\textbf{vs DE}} & \multicolumn{3}{c}{\textbf{GM Run Time (s)}} \\
         & $\mathbf{R_{\text{time}}}$ & \textbf{\textit{p}-value} & $\mathbf{R_{\text{time}}}$ & \textbf{\textit{p}-value} & \textbf{QUASAR} & \textbf{L-SHADE} & \textbf{DE} \\
        \midrule \midrule
        $\mathbf{D=10}$ & 4.21 & $5\times 10^{-144}$ & 1.26 & $7.3\times 10^{-141}$ & 16.4 & 69.1 & 20.7 \\
        $\mathbf{D=30}$ & 3.17 & $5\times 10^{-144}$ & 1.17 & $1.1\times 10^{-139}$ & 25.0 & 79.4 & 29.2 \\
        $\mathbf{D=50}$ & 2.66 & $5\times 10^{-144}$ & 1.14 & $1.6\times 10^{-137}$ & 33.5 & 89.0 & 38.1 \\
        $\mathbf{D=100}$ & 2.03 & $5\times 10^{-144}$ & 1.09 & $3\times 10^{-137}$ & 55.4 & 112.4 & 60.5 \\
        \midrule
        $\mathbf{N=100}$ & 1.37 & $1\times 10^{-142}$ & 1.05 & $1\times 10^{-84}$ & 5.8 & 7.9 & 6.0 \\
        $\mathbf{N=250}$ & 1.47 & $1\times 10^{-142}$ & 0.84 & $5\times 10^{-58}$ & 13.7 & 20.1 & 11.4 \\
        $\mathbf{N=500}$ & 1.65 & $1\times 10^{-142}$ & 1.06 & $2\times 10^{-123}$ & 28.7 & 47.2 & 30.5 \\
        $\mathbf{N=1000}$ & 2.00 & $1\times 10^{-142}$ & 1.08 & $9\times 10^{-121}$ & 57.3 & 114.4 & 62.0 \\
        $\mathbf{N=5000}$ & 5.96 & $3\times 10^{-49}$ & 1.29 & $3\times 10^{-48}$ & 292.0 & 1741.7 & 376.1 \\
        \bottomrule
    \end{tabular}
\end{table*}

\section{Ease of Use and Implementation}
\subsection{Ease of Use}
QUASAR operates most effectively with its default parameters, bypassing the hyperparameter tuning often required by other metaheuristics. The algorithm relies on its quasi-adaptive mechanisms to dynamically balance exploration and exploitation. 

Typical use cases require only the objective function, search bounds, and the target number of generations ($g_{\max}$). For additional control, the \textit{entangle rate} can be tuned between $[0,1]$, though the default value of 0.33 was most robust across the CEC benchmarks.

While not used in this study, the optimizer naturally handles vectorized objective functions by evolving the full population at once. This can lead to orders of magnitude faster run times, depending on population size.

\subsection{Implementation}
QUASAR is implemented as part of a streamlined open-source Python package, designed for ease of use and adoption into existing workflows. The syntax follows the structure of \texttt{SciPy}'s \texttt{optimize} module, making it nearly one-to-one interchangeable with its Differential Evolution function.

The source code requires minimal dependencies, relying only on \texttt{NumPy} and \texttt{SciPy}. The open-source package containing QUASAR, \texttt{hdim\_opt}, is available at \texttt{https://github.com/jgsoltes/hdim-opt,} which includes the experimental notebooks for full reproducibility.

\section{Conclusion}
This study introduced Quasi-Adaptive Search with Asymptotic Reinitialization (QUASAR), a quantum-inspired evolutionary algorithm designed to efficiently explore and exploit complex high-dimensional spaces. 

Combining probabilistic mutation, rank-based crossover, and asymptotically decaying covariance reinitializations, QUASAR’s highly stochastic architecture is exceptionally suited for discontinuous and non-parametric optimization. Its reliance on default parameters and streamlined implementation allows for seamless integration within existing workflows.

Experimental trials across the CEC2017 benchmark suite confirm that QUASAR delivers exceptional performance, achieving the highest rank via the Friedman test compared to Differential Evolution (DE) and L-SHADE, with respective overall geometric mean error improvements of $3.85\times$ and $2.07\times$. 

QUASAR also proves computationally efficient, with overall geometric mean run time improvements of $1.40\times$ and $5.16\times$ over DE and L-SHADE. Analysis by dimension and population size shows linear improvement scaling across all metrics.

The empirical evidence validates QUASAR as a powerful evolutionary algorithm, ready for the next generation of high-dimensional optimization problems. 

\section{Future Works}
QUASAR's quasi-adaptive mechanisms may be improved via strategies found in other modern DE variants.

\subsection{Enhancements}
\subsubsection{Mutations:}
The mutation factors are sampled from fixed Gaussian distributions; dynamic scaling based on convergence may significantly help to maintain a healthy balance of exploration and exploitation. The strategies may benefit from different permutations of equation terms, or non-probablistic selction methods.

\subsubsection{Crossover:}
Crossover rate is currently defined by fitness rank, and may benefit from a probability-based approach similar to the other QUASAR components, via soft-max or similar methods. Many techniques from existing evolutionary frameworks may prove beneficial to the model.

\subsubsection{Reinitialization:}
The current reinitialization covariance matrix ($\mathbf{\Sigma}$) treats all top solutions equally; future work will investigate weighted solutions similar to those found in CMA-ES. New methods for defining the probability of occurrence may be explored, such as sinusoidal or sigmoid activation behavior. 

Additionally, the reinitialization is currently bypassed for problems where $N<D$. Testing other redistribution methods, such as uniform quasi-Monte Carlo sequences, may be beneficial for these cases.

\subsubsection{Verification:} 
While QUASAR's highly stochastic framework makes it a practical choice for non-differentiable and non-parametric problems, the algorithms’s performance would be further validated by rigorously evaluating its performance on constrained or multi-objective problems.

\section*{Acknowledgements}
The author would like to thank Dr. Kellen Sorauf, Dr. Ksenia Polson, Dr. Mike Busch, and Monet Morris for their invaluable feedback and support throughout this study. It would not have been possible without the insights gained from each of our discussions.

Dr. Sorauf's recommendation to include L-SHADE in the experiment significantly improved the rigor of the study; Dr. Polson's lectures laid the framework for the statistical analysis; Dr. Busch pushed me to pursue the project to begin with. Monet remains a source of inspiration, as always.

This research was conducted independently and was not part of a funded project.

\bibliographystyle{elsarticle-num}
\bibliography{myreferences.bib}

\end{document}